\documentclass{amsart}
\usepackage{amssymb, epsf, amscd}
\allowdisplaybreaks
\swapnumbers
\newtheorem{theorem}[subsection]{Theorem}               
               
\newtheorem{acknowledgements}[subsection]{Acknowledgements}

\newtheorem{question}[subsection]{Question}         
\theoremstyle{definition}                      
\newtheorem{definition}[subsection]{Definition}

\theoremstyle{remark}

\numberwithin{equation}{section}
\def\Arc[#1]{
\ifcase#1
\qbezier[25](0.966,-0.259)(1.04,0)(0.966,0.259)
\or
\qbezier[25](0.966,0.259)(0.897,0.518)(0.707,0.707)
\or
\qbezier[25](0.707,0.707)(0.518,0.897)(0.259,0.966)
\or
\qbezier[25](0.259,0.966)(0,1.04)(-0.259,0.966)
\or
\qbezier[25](-0.259,0.966)(-0.518,0.897)(-0.707,0.707)
\or
\qbezier[25](-0.707,0.707)(-0.897,0.518)(-0.966,0.259)
\or
\qbezier[25](-0.966,0.259)(-1.04,0)(-0.966,-0.259)
\or
\qbezier[25](-0.966,-0.259)(-0.897,-0.518)(-0.707,-0.707)
\or
\qbezier[25](-0.707,-0.707)(-0.518,-0.897)(-0.259,-0.966)
\or
\qbezier[25](-0.259,-0.966)(0,-1.04)(0.259,-0.966)
\or
\qbezier[25](0.259,-0.966)(0.518,-0.897)(0.707,-0.707)
\or
\qbezier[25](0.707,-0.707)(0.897,-0.518)(0.966,-0.259)
\fi}
\def\Sterm{
\DottedCircle
\Arc[8]
\Arc[9]
\Arc[10]
\qbezier[10](0,-1)(0,-0.75)(0,-0.5)
\qbezier[10](0,-0.5)(0.25,-0.25)(0.5,0)
\put(-0.5,-0.866){\vector(1,-1){0.01}}
\qbezier[10](0,-0.5)(-0.25,-0.25)(-0.5,0)}
\def\Tterm{
\DottedCircle
\Arc[8]
\Arc[9]
\Arc[10]
\put(-0.5,-0.866){\vector(1,-1){0.01}}
\qbezier[18](-0.5,0)(-0.1,-0.5)(-0.1,-0.97)
\qbezier[18](+0.5,0)(+0.1,-0.5)(+0.1,-0.97)
}
\def\Uterm{
\DottedCircle
\Arc[8]
\Arc[9]
\Arc[10]
\put(-0.5,-0.866){\vector(1,-1){0.01}}
\qbezier[18](-0.5,0)(0.1,-0.5)(0.1,-0.97)
\qbezier[18](+0.5,0)(-0.1,-0.5)(-.1,-0.97)
}
\def\DottedCircle{
\qbezier[4](0.966,-0.259)(1.04,0)(0.966,0.259)
\qbezier[4](0.966,0.259)(0.897,0.518)(0.707,0.707)
\qbezier[4](0.707,0.707)(0.518,0.897)(0.259,0.966)
\qbezier[4](0.259,0.966)(0,1.04)(-0.259,0.966)
\qbezier[4](-0.259,0.966)(-0.518,0.897)(-0.707,0.707)
\qbezier[4](-0.707,0.707)(-0.897,0.518)(-0.966,0.259)
\qbezier[4](-0.966,0.259)(-1.04,0)(-0.966,-0.259)
\qbezier[4](-0.966,-0.259)(-0.897,-0.518)(-0.707,-0.707)
\qbezier[4](-0.707,-0.707)(-0.518,-0.897)(-0.259,-0.966)
\qbezier[4](-0.259,-0.966)(0,-1.04)(0.259,-0.966)
\qbezier[4](0.259,-0.966)(0.518,-0.897)(0.707,-0.707)
\qbezier[4](0.707,-0.707)(0.897,-0.518)(0.966,-0.259)
}

\newcommand{\firef}[1]{Figure~{\rm\ref{#1}}}
\newcommand{\thref}[1]{Theorem~\ref{#1}}

\newcommand{\fig}[1]
        {\raisebox{-0.5\height}%
                 {\epsfbox{#1}}%
        }

\newcommand\lbb[1]{\label{#1} 
                   }                                    
\def\Zset{\mathbb{Z}}       
\def\Qset{\mathbb{Q}}       
       
\def\Ct{{C}}
\def\CK{\mbox{K}}
\def\be{\beta}

\def\be{\begin{equation}}
\def\fe{\end{equation}}
\def\KI{Z}
\def\wh{\omega}
\def\CA{A}
\def\Picture#1{
\begin{picture}(2,1)(-1,-0.167)
#1
\end{picture}
}

\begin{document}

\title{Covering spaces over claspered knots}
\author{Andrew Kricker}
\address{Department of Mathematical and Computing Science, Tokyo Institute of Technology}
\email{kricker@is.titech.ac.jp}
\date{First version: December , 1998. Updated version: February , 2000}

\begin{abstract}

In this note we reconsider a familiar result in Vassiliev knot
theory - that the coefficients of the Alexander-Conway polynomial determine 
the top row of the Kontsevich integral - from the point of view of Kazuo
Habiro's clasper theory. We observe that in this setting the 
calculation reflects the topology of the universal cyclic covering space
of a claspered knot's complement.

\end{abstract}
\maketitle


\section{Introduction}\lbb{sintro}

In recent times the focus of the finite type theory of knots has shifted
to questions of what of the topology of knots finite type invariants
actually detect \cite{Hab1,Hab2,KL,N,S}. In this respect, it seems that Kazuo
Habiro's clasper theory will provide a powerful tool. This note will show
how clasper theory reveals topological content in a previously 
combinatorial understanding: that the Alexander-Conway polynomial determines
the top row of the Kontsevich integral \cite{C,K,KSA,V}, following the 
basic framework of \cite{BNG}. The original conjecture in whose shadow
this work lies is due to Melvin and Morton \cite{MM}. 

We note that certain aspects of Habiro's clasper theory have
been independently considered by Goussarov. In particular see the
paper \cite{G} and the web page \cite{GW}. 

Take the following representation of the Alexander-Conway polynomial (as a
formal power series in the variable $h$):
\begin{eqnarray*}
\Ct_{K^+}(h) - \Ct_{K^-}(h) & = & (e^{\frac{h}{2}}-e^{-\frac{h}{2}})\Ct_{K^s}(h), \\
\Ct_U(h)& = & 1.
\end{eqnarray*}
In the usual way $K^+,K^-$ and $K^s$ are three knots differing in a
ball as an overcrossing, an undercrossing, and an orientation-preserving 
smoothing of the crossing, and $U$ is the unknot.
Take the coefficients of this series
$\Ct_K(h)=\sum_{n=0}^{\infty} {c}_{2n}(K) h^{2n}$.
The coefficient ${c}_{2n+1}$ vanishes and the coeffiient ${c}_{2n}$ is finite
type of order $2n$ \cite{BN1}. Let the coefficients determine rational knot invariants $d_{2n}(K)$:
\begin{equation}
\Ct_{K}(h) = \exp{(-2\sum_{n=1}^{\infty} d_{2n}(K) h^{2n})}.
\end{equation}

In these terms the top row of the Kontsevich integral is determined:
\begin{equation}
\KI(K) = \exp( \frac{f(K)}{2}\theta + \sum_{n=1}^{\infty} d_{2n}(K) \wh_{2n} + X),
\end{equation}
where $\KI$ is the framed Kontsevich integral \cite{LM} with the normalisation 
employed in \cite{BN2} (so that $\KI(U)=1$),
where $\theta$ is the web diagram with a single chord, where
$\wh_{2n}$ is the ``wheel with $2n$ spokes'',\vspace{0.25cm}
\newpage
\begin{equation}
\wh_2 =\ \ \ \
\Picture{
\thicklines
\put(0,0){\circle{4}}
\DottedCircle
\qbezier[6](1.95,0)(1.5,0)(1,0)
\qbezier[6](-2,0)(-1.5,0)(-1,0)
\put(0,2){\vector(-1,0){0.01}}}\ \ \ \ \ ,\ 
\wh_4 =\ \ \ \
\Picture{
\thicklines
\put(0,0){\circle{4}}
\DottedCircle
\qbezier[6](1.95,0)(1.5,0)(1,0)
\qbezier[6](-2,0)(-1.5,0)(-1,0)
\qbezier[6](0,1)(0,1.5)(0,2)
\qbezier[6](0,-1)(0,-1.5)(0,-2)
\put(0,2){\vector(-1,0){0.01}}}\ \ \ \ \ ,\ 
\wh_6 =\ \ \ \
\Picture{
\thicklines
\put(0,0){\circle{4}}
\DottedCircle
\qbezier[6](1.95,0)(1.5,0)(1,0)
\qbezier[6](1,1.73)(0.75,1.3)(0.5,0.87)
\qbezier[6](-1,1.73)(-0.75,1.3)(-0.5,0.87)
\qbezier[6](1,-1.73)(0.75,-1.3)(0.5,-0.87)
\qbezier[6](-1,-1.73)(-0.75,-1.3)(-0.5,-0.87)
\qbezier[6](-2,0)(-1.5,0)(-1,0)
\put(0,2){\vector(-1,0){0.01}}}\ \ \ \ \ ,\ \ldots
\lbb{wheels}
\end{equation}
\[\]
and where $X$ is a series of diagrams with connected dashed graphs of negative
Euler characteristic. (The ``Chinese character diagram'' wheels first 
appeared in \cite{CV}. The papers \cite{GH} and \cite{W} contain this 
formula. )

This understanding follows from the theory of the Kontsevich integral, 
observations concerning the structure of the primitive subspace
of the chord diagram algebra \cite{CV}, and the following theorem, where 
the weight system of ${c}_{2n}$ is denoted $W({c}_{2n})$. 
Pieces of the following theorem appear in many papers \cite{BNG,C,K,KSA,V},
most significantly \cite{BNG}; 
\cite{C} and \cite{V} contain clear presentations
of the final product:

\begin{theorem}\lbb{it}
\begin{eqnarray}
& & W({c}_{2n})( \KI_{2n}(K)) = {c}_{2n}(K), \lbb{FA}\\
& & W({c}_{m_1+m_2})(D_1\#D_2) = W({c}_{m_1})(D_1)W({c}_{m_2})(D_2), \lbb{FB}\\
& & W({c}_{2n})(w_{2n})  = -2,\ \ \ \label{FC}\\
& & W({c}_{2n})(D)  =  0\ \ \ 
\begin{array}{l}\mbox{if the dashed graph of $D$} \\
\mbox{has -ve Euler characteristic,}
\end{array} \lbb{vanish}\lbb{FD}
\end{eqnarray}
where $\KI_{2n}$ is the degree $2n$ term of $\KI$, $D_i\in \CA_{m_i}(S^1)$ and $D$ is a diagram in $\mbox{Prim}(\CA_{2n}(S^1))$.
\end{theorem}

A major virtue of clasper theory is that it constructs directly from some 
primitive web diagram $D$, a surgery presentation of a knot $K$ which has the
property that 
$\KI(K) = 1 + D\ +\ $higher order terms. Weight systems of invariants can thus
be calculated on primitive web diagrams by evaluating the difference between 
the value the invariant takes on a single knot claspered from the unknot
according to the diagram $D$  and the value it takes on the unknot
(as opposed to the usual double sum that comes from STU resolutions
of trivalent vertices of the dashed graph followed by desingularisations of 
double points).

In this note the Alexander-Conway invariants of knots obtained by claspering
the unknot
will be calculated by explicitly constructing surgery descriptions of the
universal cyclic covers of their complements. The knots associated to wheels are familiar as twisted
variants of Ng's wheels \cite{N} (and in fact, the approach presented in this
work owes a debt to \cite{N}). In this framework the vanishing theorem (Equation~{\rm\ref{FD}}) recieves a topological explanation. This approach provokes
some questions. 

Generalisations of some of these
observations in the direction of knot invariants which factor through 
Seifert matrices
are contained in recent work due to Murakami and Ohtsuki \cite{MO}. There are close relations with some extensive 
recent work
due to Habiro, Kanenobu and Shima concerning finite-type invariants of ribbon 2-knots \cite{HKS,HS}. 

\begin{acknowledgements}
The author is supported by a Japan Society for the Promotion of Science
postdoctoral fellowship. It is a pleasure to acknowledge discussions with
Kazuo Habiro and Tomotada Ohtsuki, and a reading of an early version by
Simon Willerton.
\end{acknowledgements}
\section{Finite type invariants and clasper theory}\lbb{clathpers}

The finite type theory of knots is based on a filtration of the $\Qset$-vector
space of finite, $\Qset$-linear combinations of knots, ${\CK}$. Here this 
filtration will be introduced by means of {\it claspers}. 
Clasper theory is due to Kazuo Habiro \cite{Hab1,Hab2} (extensively
generalising precursors due to Matveev \cite{Ma} and
Murakami--Nakanishi \cite{MN}).
We note that certain aspects of this theory were independently
considered by Goussarov \cite{G,GW}.
A clasper is an embedding
of a band-summed pair of annuli in the complement of the knot.

\begin{equation}
\setlength{\unitlength}{30pt}
\Picture{
\qbezier(1,-0.5)(1,-0.8)(1,-1.4)
\qbezier(1,1)(1,0.5)(1,-0.3)
\qbezier(-1,1)(-1,0.5)(-1,0.1)
\qbezier(-1,-1.4)(-1,-1.0)(-1,-0.1)
\put(-1,0.8){\vector(0,-1){0.01}}
\put(1,0.8){\vector(0,1){0.01}}
\qbezier(1.1,0.1)(1.4,0.1)(1.4,-0.2)
\qbezier(1.1,-0.5)(1.4,-0.5)(1.4,-0.2)
\qbezier(1.1,-0.1)(1.2,-0.1)(1.2,-0.2)
\qbezier(1.1,-0.3)(1.2,-0.3)(1.2,-0.2)
\qbezier(1.1,-0.5)(0.9,-0.5)(0.9,-0.5)
\qbezier(1.1,-0.3)(0.9,-0.3)(0.9,-0.3)
\qbezier(0.9,-0.5)(0.6,-0.5)(0.6,-0.2)
\qbezier(0.9,-0.3)(0.8,-0.3)(0.8,-0.2)
\qbezier(0.9,-0.1)(0.8,-0.1)(0.8,-0.2)
\qbezier(0.9,0.1)(0.6,0.1)(0.6,-0.2)
\qbezier(-1.1,-0.5)(-1.4,-0.5)(-1.4,-0.2)
\qbezier(-1.1,0.1)(-1.4,0.1)(-1.4,-0.2)
\qbezier(-1.1,-0.3)(-1.2,-0.3)(-1.2,-0.2)
\qbezier(-1.1,-0.1)(-1.2,-0.1)(-1.2,-0.2)
\qbezier(-1.1,0.1)(-0.9,0.1)(-0.9,0.1)
\qbezier(-1.1,-0.1)(-0.9,-0.1)(-0.9,-0.1)
\qbezier(-0.9,0.1)(-0.6,0.1)(-0.6,-0.2)
\qbezier(-0.9,-0.1)(-0.8,-0.1)(-0.8,-0.2)
\qbezier(-0.9,-0.3)(-0.8,-0.3)(-0.8,-0.2)
\qbezier(-0.9,-0.5)(-0.6,-0.5)(-0.6,-0.2)
\qbezier(0.6,-0.1)(-0,-0.1)(-0.6,-0.1)
\qbezier(0.6,-0.3)(0,-0.3)(-0.6,-0.3)
}
\ \ \ \ \ \ \rightarrow\ \ \ \ \ 
\Picture{
\qbezier(1,-0.5)(1,-0.8)(1,-1.4)
\qbezier(1,1)(1,0.5)(1,-0.3)
\qbezier(-1,1)(-1,0.5)(-1,0.1)
\qbezier(-1,-1.4)(-1,-1.0)(-1,-0.1)
\qbezier(-1.1,0)(-1,0)(-0.1,0)
\qbezier(-0.1,0)(0.1,0)(0.1,-0.2)
\qbezier(1.1,-0.4)(0.1,-0.4)(0.1,-0.4)
\qbezier(0.1,-0.4)(-0.1,-0.4)(-0.1,-0.2)
\qbezier(0.1,-0.2)(0.1,-0.275)(0.05,-0.35)
\qbezier(-0.1,-0.2)(-0.1,-0.125)(-0.05,-0.05)
\qbezier(0.1,0)(0.1,0)(0.9,0)
\qbezier(-0.1,-0.4)(-0.1,-0.4)(-0.9,-0.4)
\qbezier(1.1,0)(1.3,0)(1.3,-0.2)
\qbezier(1.1,-0.4)(1.3,-0.4)(1.3,-0.2)
\qbezier(-1.1,0)(-1.3,0)(-1.3,-0.2)
\qbezier(-1.3,-0.2)(-1.3,-0.4)(-1.1,-0.4)
\put(0.4,0.15){\mbox{0}}
\put(-0.6,0.15){\mbox{0}}
\put(-1,0.8){\vector(0,-1){0.01}}
\put(1,0.8){\vector(0,1){0.01}}
}
\ \ \ \ \rightarrow\ \ \ \ \
\Picture{
\qbezier(1,1)(1,0.5)(1,-0.3)
\qbezier(-1,1)(-1,0.5)(-1,0.4)
\qbezier(-1,0.4)(-1,0)(0.9,0)
\qbezier(1.1,0)(1.3,0)(1.3,-0.2)
\qbezier(1.3,-0.2)(1.3,-0.4)(1,-0.4)
\qbezier(1,-0.4)(-1,-0.4)(-1,-0.8)
\qbezier(-1,-0.8)(-1,-0.8)(-1,-1.4)
\qbezier(1,-0.5)(1,-0.8)(1,-1.4)
\put(-1,0.8){\vector(0,-1){0.01}}
\put(1,0.8){\vector(0,1){0.01}}}\ \ \ \ \ .
\lbb{clathp}
\end{equation}\vspace{1cm}

A clasper carries placement information for a surgery link. In the
case of Figure \ref{clathp} the surgery link is a $0$-framed Hopf link, surgery on which
effects a crossing change on the knot. Of course, two points on the knot may
be ``claspered'' like this in many different ways corresponding to different
homotopy classes of the band. Habiro has called a clasper of this sort a 
``quantum chord''. In this work we will restrict ourselves to clasper decorations of knots such that each annulus is untwisted, and such that there is a
diagram of the decoration (in the obvious sense) with each annulus linking the
knot once as in Figure \ref{clathp} (though we will soon expand this
class to let annuli meet in something like a trivalent vertex) and with 
each band untwisted (immersed in the 
plane of the projection). This restriction essentially rules out 
half-twists of bands and is imposed here mainly to deal with sign issues
in the isomorphisms presented below.

Decorate a knot $K$ with $n$ disjoint claspers $C=\{C_i\}$, as above.
Write $K^C$ for 
the knot obtained by surgery on such a set of claspers. Denote the following
sum
\begin{equation*}
[K,C]=\sum_{D\subset C}(-1)^{\# D}K^{D} \in \CK.
\end{equation*}
\begin{definition}
${\CK}_n$ is the subspace of $\CK$ spanned by vectors $[K,C]$ where 
$\#C=n$.
\end{definition}

The graded piece $\frac{\CK_n}{\CK_{n+1}}$ is isomorphic to
$\CA'_n(S^1)$, the space of web diagrams of degree $n$ \cite{BN2}. $\CA'_n(S^1)$
is the
$\Qset$-vector space generated by web diagrams quotiented by the subspace
described by STU and 1T relations.
Web diagrams are graphs with univalent or trivalent
vertices (such as the wheels of Equation \ref{wheels}), such that each trivalent vertex is equipped with a cyclic ordering
of its incident edges and the univalent vertices are labelled (up to
orientation-preserving diffeomorphism) 
with disjoint points on an oriented copy of $S^1$. The
STU relations are
\setlength{\unitlength}{20pt}
\begin{eqnarray}
\Picture{
\Sterm} &  = &
\Picture{\Tterm}\ \ -\ \ \Picture{\Uterm}\ \ \\
& & \nonumber
\end{eqnarray}
\setlength{\unitlength}{10pt}
\hspace{-0.25cm}
where the diagrams of the three terms differ only in the parts pictured.
The 1T relations are spanned by diagrams with isolated chords. 

Call this isomorphism $\phi_n:\CA'_n(S^1)\rightarrow \frac{\CK_n}{\CK_{n+1}}$.
We will introduce Habiro's construction of this isomorphism in the next
section, using clasper graphs.

\subsection{Clasper graphs}\lbb{graphs}
Clasper theory has provided a tool to describe $\phi_n$ directly on web 
diagrams using knots decorated with {\it clasper graphs}. 
In a clasper graph, we extend the class of clasper decorations to allow
annuli to meet in a set of Borromean rings, as below (retaining the
requirement that there be a diagram such that no bands have half-twists):
\begin{equation}
\epsfxsize = 4cm
\fig{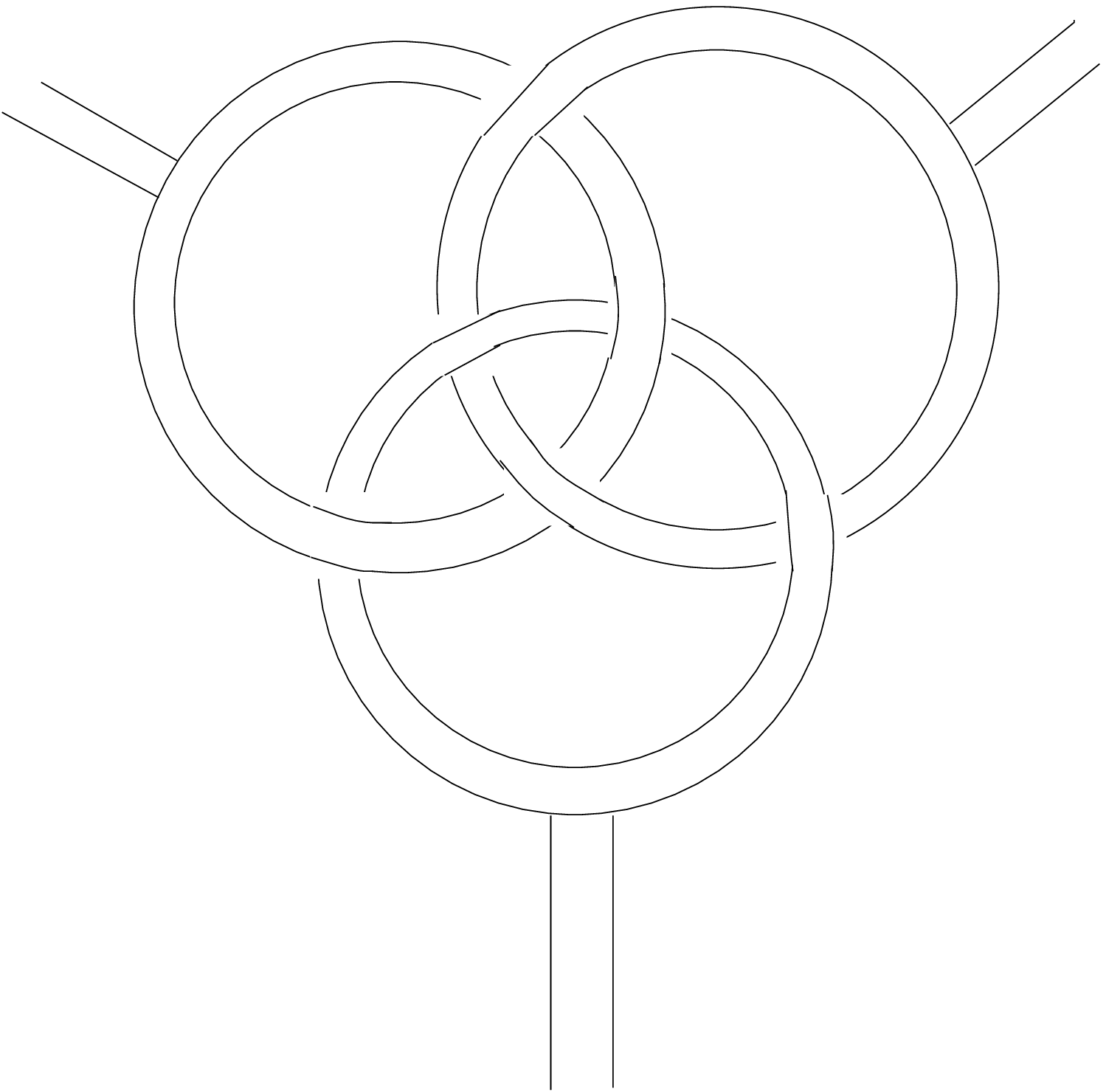} \ \ \ \ \ \ \ \ .
\end{equation}

The annuli which link the knot are called the {\it leaves} of the 
clasper graph.
Clasper theory covers a much more extensive class of objects than the class
employed in our calculation, and we enthusiastically direct the reader to 
\cite{Hab1,Hab2}.

Such a graph as we have described here has an underlying web diagram obtained
by replacing the Borromean rings with trivalent vertices (oriented according
to the diagram in which bands have no half-twists) and leaves with 
univalent vertices appropriately ordered on the $S^1$, though of
course many different clasper graphs will have the same underlying diagram.
The degree of a clasper graph is the degree of the underlying diagram (so is
half the sum of the numbers of leaves and sets of 
Borromean rings in the graph).

Decorate a knot $K$ with some 
degree $n$
clasper graph $G=\{G_1,\ldots,G_k\}$ having
underlying web diagram $D$, where the clasper 
subgraphs $\{G_1,\ldots,G_k\}$ correspond
to the connected components of the dashed graph of $D$.
\begin{theorem}[\cite{Hab2}]
\begin{equation}
\phi_n(D) = (-1)^{l+\mbox{deg}(G)}\sum_{S\subset G} (-1)^{\#G-\#S}K^{S},
\end{equation}
where $l$ is the number of leaves of $G$, the sum is over all $(2^k)$
subsets of the set
of clasper subgraphs corresponding to the connected components of the dashed graph of the underlying diagram of $G$, and $\# S$ is the cardinality of $S$.
\end{theorem}
Observe that if the underlying web diagram has connected dashed graph
(so that it is primitive) 
then its image in $\frac{\CK_n}{\CK_{n+1}}$ is a difference
$(-1)^{l+n}(K^G-K)$. 

\subsection{Remark}\lbb{technical}
There is a technical observation from clasper theory that will be required
in the sequel. Consider a knot obtained from some other knot by surgery on 
a clasper graph that is of the sort described above,
except that one clasper has a free end, as below:
\begin{equation}
\epsfxsize = 3cm
\fig{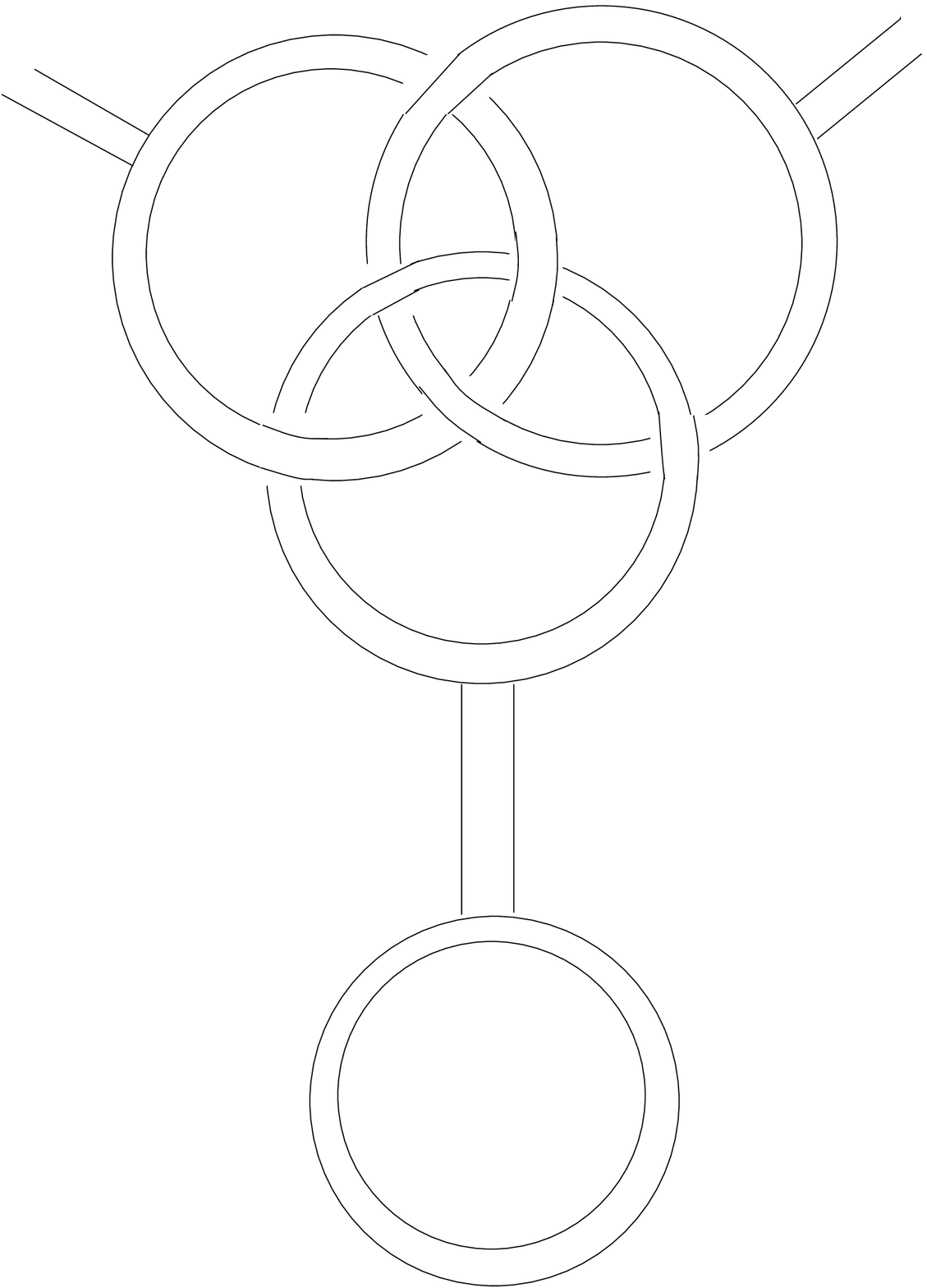}
\ \ \rightarrow\ \ \ \
\epsfxsize = 3cm 
\raisebox{1cm}{\fig{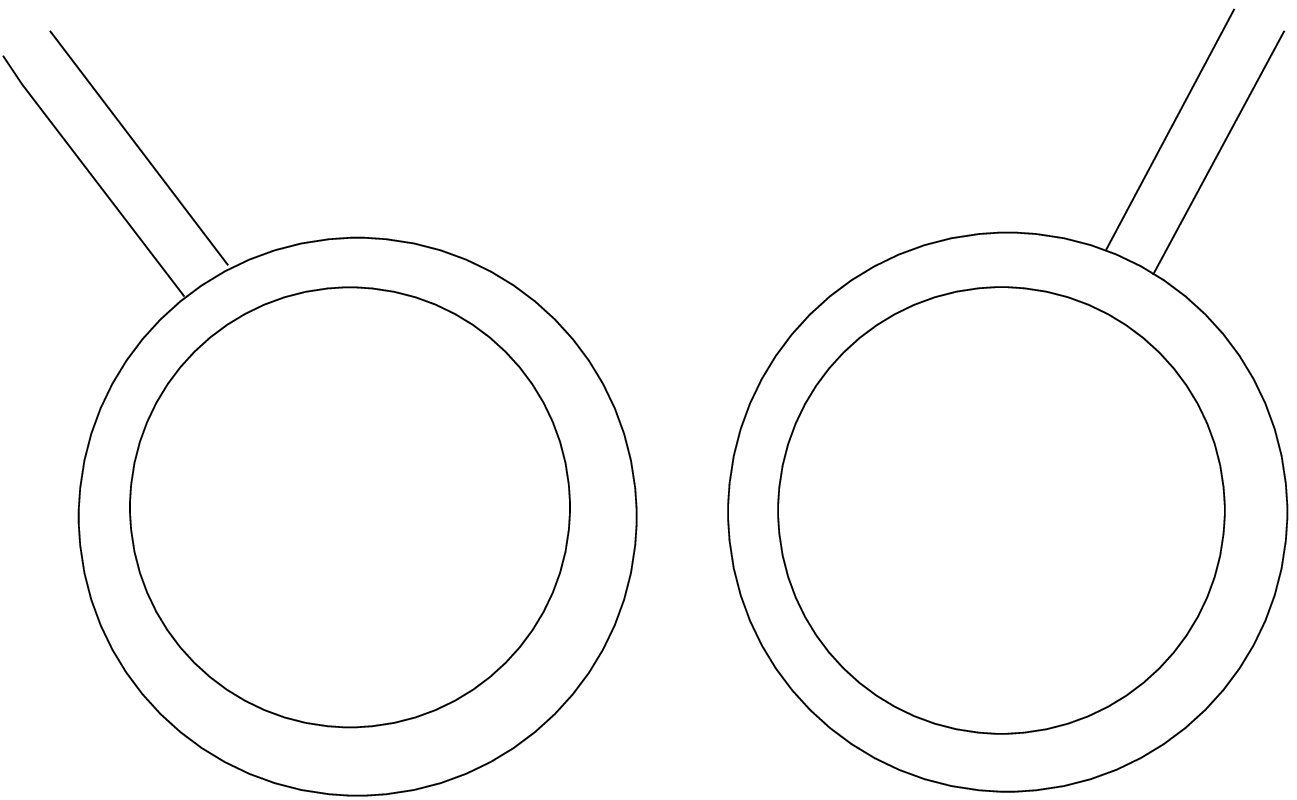}}\ \ .
\end{equation}

It is clear that a single pair of surgeries on the components associated with
the clasper with a free end will leave the graph as shown. Continuing, it is obvious that 
a component with a free end may be removed altogether without affecting the
resultant knot. (This is related to how one sees that surgery on a clasper 
graph leaves the ambient manifold invariant.)

\subsection{Weight systems}

The weight system of a degree $n$ finite type invariant $V$ (so that $V(\CK_{n+1})=0$) is the composition $W(V)=V\circ \phi_n$.

In this note it will be more convenient to use the space
of web diagrams quotiented by only the STU relations: it will 
be denoted here $\CA_n(S^1)$. 
The clasper graph construction associates to any knot invariant of finite
type of degree $n$ a dual vector to $\CA_n(S^1)$. This dual vector satisfies
1T relations and thus projects to the corresponding weight system.  
The weight system stated in the main theorem is precisely the appropriate
dual vector to $\CA_n(S^1)$ that projects to the weight system for 
$c_{2n}$. This presentation is required here as 
$\CA_n(S^1)$ is the range of the version of the
Kontsevich integral referred to in the introduction.

\section{Covering spaces of claspered knots}

\subsection{The Alexander polynomial}

A short review is in order. A standard reference is \cite{Rolf}.
Take $\widetilde{S^3\setminus K}$, the covering space
of $S^3\setminus K$ corresponding to the abelianisation of 
$\pi = \pi_1(S^3\setminus K)$ (that is, to the short exact sequence
$0\rightarrow [\pi,\pi]\rightarrow \pi\rightarrow \Zset\rightarrow 0$),
with group of deck transformations $\Zset$ generated by some $T$. 
This is called the universal cyclic covering space. Defining
a $\Zset[t,t^{-1}]$ action on $H_1(\widetilde{S^3\setminus K},\Zset)$ by $\sum_i a_i t^{n_i}.h
=\sum_i a_i (T_*)^{n_i}(h)$, $H_1(\widetilde{S^3\setminus K},\Zset)$ becomes a
$\Zset[t,t^{-1}]$-module. The Alexander polynomial of $K$ (up to multiplication by 
units) is a generator of the order ideal of $H_1(\widetilde{S^3\setminus K},\Zset)$.

\subsection{Complete diagrams}

Let us call a diagram {\it complete} if the trivalent vertices of its dashed graph have at most one incident edge whose other end is a univalent vertex (this
extends a term due to Ng). A primitive web diagram is equal to a sum of 
complete primitive web diagrams of lesser or equal Euler characteristic. Thus 
we may restrict out attention to complete diagrams in the proof of \thref{it}.

\subsection{Universal cyclic covering spaces}\lbb{uccs}

Take a complete primitive web diagram $D\in \CA_n(S^1)$. Decorate the
unknot with a clasper graph $G$ with underlying diagram $D$: so that $(-1)^{l+n} (U^G-U)$ represents
$\phi_n(D)$ in $\frac{\CK_n}{\CK_{n+1}}$. 
Choose $G$ so that if the leaves are removed then
the resulting graph is unlinked from the knot (this can always be done, and is
done here essentially to simplify the choice of Seifert surface in the exposition
below).
$W_n({c}_n)(D)$ is equal to $(-1)^{l+n}{c}_n(U^G)$.
We calculate this by building a surgery description of the universal
cyclic cover of $S^3\setminus U^G$.

Observe that near a leaf of the clasper graph, a single surgery on the clasper
ending in the leaf brings the picture to \firef{seifert}.\vspace{2cm}

\begin{equation}\lbb{seifert}
\setlength{\unitlength}{20pt}
\Picture{
\qbezier(-1.9,0.3)(-0.7,0)(1.3,-0.5)
\put(-1.1,0.1){\vector(4,-1){0.01}}
\qbezier[10](-0.2,-.125)(0.05,0.56)(0.3,1)
\qbezier[10](0.3,1)(0.3,1.5)(0.3,2)
\qbezier[10](0.3,1)(0.8,1)(1.3,1)}
\ \ \ \ \ \ \ \ \rightarrow\ \ \ \ \ \ \ \ \ \ \ \ \ \ \ \ \ \ \ \ \ \ \ \
\setlength{\unitlength}{30pt}
\Picture{
\qbezier(0,0)(0,1)(0,1)
\qbezier(0,0)(0,-0.3)(0.3,-0.3)
\qbezier(-0.5,0)(-0.5,1)(-0.5,1)
\qbezier(-0.5,0)(-0.5,-0.3)(-0.8,-0.3)
\qbezier(-0.5,1.2)(-0.5,1.4)(-0.5,1.4)
\qbezier(0,1.2)(0,1.4)(0,1.4)
\qbezier(-0.5,1.4)(-0.5,1.8)(-1.2,1.8)
\qbezier(0,1.4)(0,2.3)(-1.2,2.3)
\qbezier(-1.2,1.8)(-1.9,1.8)(-1.9,1.4)
\qbezier(-1.2,2.3)(-2.4,2.3)(-2.4,1.4)
\qbezier(-1.9,0)(-1.9,-0.3)(-1.6,-0.3)
\qbezier(-2.4,0)(-2.4,-0.3)(-2.7,-0.3)
\qbezier(0,0.5)(0,0.5)(3.8,-0.45)
\qbezier(-0.5,0.625)(-0.5,0.625)(-1.9,0.975)
\qbezier(-2.4,1.1)(-2.4,1.1)(-3.2,1.3)
\put(-3.6,0.8){$\Sigma$}
\qbezier(-1.9,1.4)(-1.9,1.4)(-1.9,0)
\qbezier(-2.4,1.4)(-2.4,1.4)(-2.4,0)
\qbezier(-0.6,1.5)(-0.8,1.5)(-0.8,1.25)
\qbezier(-0.8,1.25)(-0.8,1.1)(-0.5,1.1)
\qbezier(-0.5,1.1)(1.2,1.1)(1.2,1.1)
\qbezier(1.2,1.1)(1.5,1.1)(1.5,1.4)
\qbezier(1.5,1.6)(1.5,2.0)(1.5,2.0)
\qbezier(1.5,2.2)(1.5,2.4)(1.3,2.4)
\qbezier(1.3,2.4)(1.1,2.4)(1.1,2.2)
\qbezier(0.8,1.5)(1.1,1.5)(1.1,1.8)
\qbezier(1.1,1.8)(1.1,2.2)(1.1,2.2)
\qbezier(0,1.5)(0,1.5)(0.8,1.5)
\qbezier(1.3,2.45)(1.3,2.45)(1.3,3)
\qbezier(1.3,2.35)(1.3,2.1)(1.5,2.1)
\qbezier(1.5,2.1)(1.7,2.1)(1.7,2.4)
\qbezier(1.7,2.4)(1.7,3)(1.7,3)
\qbezier(1.5,1.5)(1.1,1.5)(1.1,1.15)
\qbezier(1.1,1.05)(1.1,1.05)(1.1,0.28)
\qbezier(1.5,1.5)(1.5,1.5)(2.5,1.5)
\qbezier(2.5,1.5)(2.7,1.5)(2.7,1.3)
\qbezier(2.7,1.3)(2.7,1.1)(2.5,1.1)
\qbezier(2.3,1.1)(2.3,1.1)(1.7,1.1)
\qbezier(1.7,1.1)(1.5,1.1)(1.5,0.9)
\qbezier(1.5,0.9)(1.5,0.17)(1.5,0.17)
\qbezier[10](1.1,0.1)(1.1,-0.6)(1.1,-1.2)
\qbezier[10](1.5,0.0)(1.5,-0.6)(1.5,-1.2)
\qbezier[6](1.1,-1.2)(1.1,-1.4)(0.9,-1.4)
\qbezier[10](1.5,-1.2)(1.5,-1.8)(0.9,-1.8)
\qbezier[6](-0.25,-1.2)(-0.25,-1.4)(-0.05,-1.4)
\qbezier[10](-2.15,-1.2)(-2.15,-1.8)(-1.55,-1.8)
\qbezier[16](-0.05,-1.4)(0.4,-1.4)(0.9,-1.4)
\qbezier[24](-1.55,-1.8)(-1.2,-1.8)(0.9,-1.8)
\qbezier[24](-2.15,-1.2)(-2.15,0.1)(-2.15,1.4)
\qbezier[20](-0.25,-1.2)(-0.25,0.1)(-0.25,1)
\qbezier[16](-2.15,1.4)(-2.15,2.05)(-1.2,2.05)
\qbezier[20](-1.2,2.05)(-0.25,2.05)(-0.25,1.2)
\qbezier(2.4,1.1)(2.4,1.3)(2.65,1.3)
\qbezier(2.4,1.1)(2.4,0.9)(2.65,0.9)
\qbezier(2.75,1.3)(3.5,1.3)(3.5,1.3)
\qbezier(2.65,0.9)(3.5,0.9)(3.5,0.9)
\put(0.8,2.7){$0$}
\put(3.2,1.45){$0$}
\put(-1.2,1.25){$0$}
\put(0.7,0.65){$0$}
}
\setlength{\unitlength}{10pt}
\end{equation}\vspace{1.3cm}

In \firef{seifert} a Seifert surface $\Sigma$ for the unknot which is 
disjoint from the surgery link has been chosen.
Use this surface to construct the universal cyclic cover of the
unknot in the usual way. [ In outline:
we ``split'' $S^3\setminus U$ along $\Sigma\setminus\partial\Sigma$ (which we
will denote $ $$^\circ\Sigma$) leaving
a manifold $X$ with boundary $ $$^\circ\Sigma^{\pm}$, two copies of
$ $$^\circ\Sigma$ which meet along $U$ (although $U$ is excluded from the space).
$S^3\setminus U$ is recovered from an identification $\sigma: $$^\circ\Sigma^+\rightarrow
 $$^\circ\Sigma^-$. Take $\Zset$ copies of the triple $( $$^\circ\Sigma^-,X, $$^\circ\Sigma^+)$ denoting
them $( $$^\circ\Sigma^-_{i},X_i, $$^\circ\Sigma^+_{i})$ and $\Zset$ of the appropriate maps
$ $$^\circ\Sigma^+_i\rightarrow  $$^\circ\Sigma^{-}_{i+1}$. $\widetilde{S^3\setminus U}$ is recovered
as the identification space of this system.
A rather more rigorous discussion of this procedure may
be found in \cite{Rolf}].

Denote the covering map $p:\widetilde{S^3\setminus U}
\rightarrow S^3\setminus U$. As $\Sigma$ is disjoint from the surgery link it is clear
that $L$ lifts, $p^{-1}(L)$, to a surgery presentation for
$\widetilde{S^3\setminus U^G}$.
\subsection{Proof of the vanishing theorem, Equation~{\rm\ref{vanish}}}

Take $G$ some graph clasper for $U$ with underlying diagram some chosen complete 
primitive diagram $D$ of negative Euler characteristic. The dashed graph of
such a diagram will 
have at least one trivalent vertex which has no incident edges connecting 
univalent vertices. Thus the corresponding surgery link will have at least one
set of three components linking in a set of Borromean rings in $X_0$, lifted
to all the $X_i$. But the Alexander polynomial is calculated from the homology
of $\widetilde{S^3\setminus U^G}$ and this coincides with the homology of
$\widetilde{S^3\setminus U^{G'}}$ where $G'$ is obtained from $G$ by unlinking that set of Borromean links:
\begin{equation}
\epsfxsize = 4cm
\fig{clasp3.eps}
\ \ \rightarrow\ \
\epsfxsize = 4cm
\fig{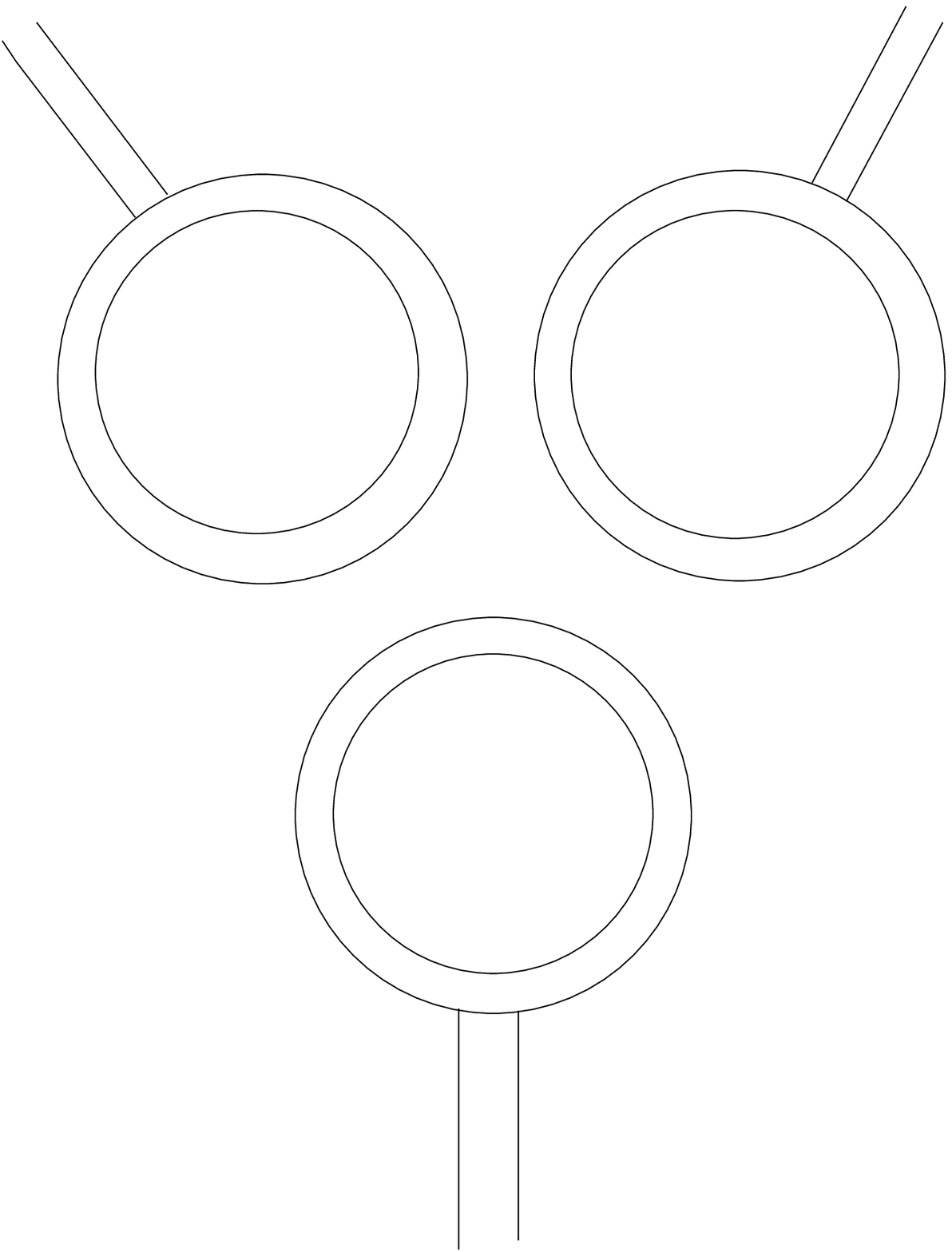}\ \ \ .
\end{equation}

Further, the observation in Remark~{\rm\ref{technical}} shows that $U^{G'}$ is in fact the unknot,
so that $U^G$ has trivial Alexander polynomial.

\subsection{Wheels}

We turn our attention to the wheels. Denote by $W_{2n}$ the graph clasper for 
the unknot with underlying graph corresponding to the 
planar embedding of the wheel $\omega_{2n}$. For example, $(U,W_{2})$ is
\begin{equation}
\epsfxsize = 4cm
{\fig{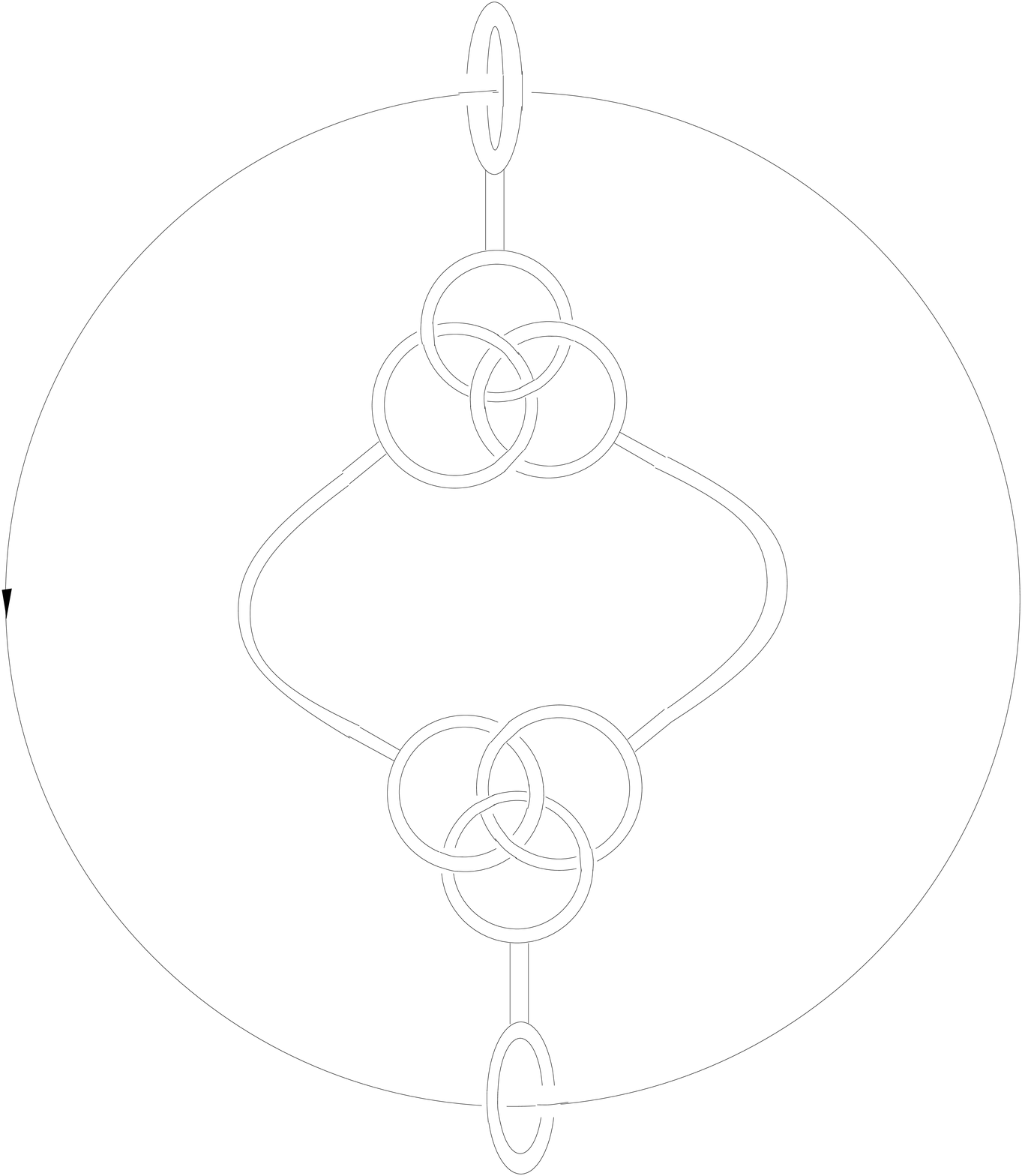}}\ \ .
\end{equation}
This sequence represents the wheel diagrams in the associated
graded vector space: $\phi_{2n}(\wh_{2n})=U^{W_{2n}}-U$. 

Number the leaves of $W_{2n}$
clockwise $1\ldots 2n$.
Around leaf $i$ of the clasper graph the picture looks as follows (where
$i+1(2n)$ means $i+1$ modulo $2n$, etc):

\setlength{\unitlength}{25pt}
\begin{eqnarray}
& &\lbb{wheelv}
\ \ \ \ \ \ \ \ \ \ \ \ \ \ \ \ \ \ \Picture{
\qbezier(0,0)(0,1)(0,1)
\qbezier(0,0)(0,-0.3)(0.3,-0.3)
\qbezier(-0.5,0)(-0.5,1)(-0.5,1)
\qbezier(-0.5,0)(-0.5,-0.3)(-0.8,-0.3)
\qbezier(-0.5,1.2)(-0.5,1.4)(-0.5,1.4)
\qbezier(0,1.2)(0,1.4)(0,1.4)
\qbezier(-0.5,1.4)(-0.5,1.8)(-1.2,1.8)
\qbezier(0,1.4)(0,2.3)(-1.2,2.3)
\qbezier(-1.2,1.8)(-1.9,1.8)(-1.9,1.4)
\qbezier(-1.2,2.3)(-2.4,2.3)(-2.4,1.4)
\qbezier(-1.9,0)(-1.9,-0.3)(-1.6,-0.3)
\qbezier(-2.4,0)(-2.4,-0.3)(-2.7,-0.3)
\qbezier(0,0.5)(0,0.5)(3,-0.25)
\qbezier(-0.5,0.625)(-0.5,0.625)(-1.9,0.975)
\qbezier(-2.4,1.1)(-2.4,1.1)(-3.2,1.3)
\qbezier(-1.9,1.4)(-1.9,1.4)(-1.9,0)
\qbezier(-2.4,1.4)(-2.4,1.4)(-2.4,0)
\qbezier(-0.6,1.5)(-0.8,1.5)(-0.8,1.25)
\qbezier(-0.8,1.25)(-0.8,1.1)(-0.5,1.1)
\qbezier(-0.5,1.1)(1.2,1.1)(1.2,1.1)
\qbezier(1.2,1.1)(1.5,1.1)(1.5,1.4)
\qbezier(1.5,1.6)(1.5,2.0)(1.5,2.0)
\qbezier(1.5,2.2)(1.5,2.4)(1.3,2.4)
\qbezier(1.3,2.4)(1.1,2.4)(1.1,2.2)
\qbezier(0.8,1.5)(1.1,1.5)(1.1,1.8)
\qbezier(1.1,1.8)(1.1,2.2)(1.1,2.2)
\qbezier(0,1.5)(0,1.5)(0.8,1.5)
\qbezier(1.3,2.45)(1.3,2.45)(1.3,3)
\qbezier(1.3,2.35)(1.3,2.1)(1.5,2.1)
\qbezier(1.5,2.1)(1.7,2.1)(1.7,2.4)
\qbezier(1.7,2.4)(1.7,3)(1.7,3)
\qbezier(1.5,1.5)(1.1,1.5)(1.1,1.15)
\qbezier(1.1,1.05)(1.1,1.05)(1.1,0.28)
\qbezier(1.5,1.5)(1.5,1.5)(2.5,1.5)
\qbezier(2.5,1.5)(2.7,1.5)(2.7,1.3)
\qbezier(2.7,1.3)(2.7,1.1)(2.5,1.1)
\qbezier(2.3,1.1)(2.3,1.1)(1.7,1.1)
\qbezier(1.7,1.1)(1.5,1.1)(1.5,0.9)
\qbezier(1.5,0.9)(1.5,0.17)(1.5,0.17)
\qbezier[10](1.1,0.1)(1.1,-0.6)(1.1,-1.2)
\qbezier[10](1.5,0.0)(1.5,-0.6)(1.5,-1.2)
\qbezier[6](1.1,-1.2)(1.1,-1.4)(0.9,-1.4)
\qbezier[10](1.5,-1.2)(1.5,-1.8)(0.9,-1.8)
\qbezier[6](-0.25,-1.2)(-0.25,-1.4)(-0.05,-1.4)
\qbezier[10](-2.15,-1.2)(-2.15,-1.8)(-1.55,-1.8)
\qbezier[16](-0.05,-1.4)(0.4,-1.4)(0.9,-1.4)
\qbezier[24](-1.55,-1.8)(-1.2,-1.8)(0.9,-1.8)
\qbezier[24](-2.15,-1.2)(-2.15,0.1)(-2.15,1.4)
\qbezier[20](-0.25,-1.2)(-0.25,0.1)(-0.25,1)
\qbezier[16](-2.15,1.4)(-2.15,2.05)(-1.2,2.05)
\qbezier[20](-1.2,2.05)(-0.25,2.05)(-0.25,1.2)
\qbezier(2.4,1.1)(2.4,1.3)(2.65,1.3)
\qbezier(2.4,1.1)(2.4,0.9)(2.65,0.9)
\qbezier(2.75,1.3)(3.5,1.3)(3.5,1.3)
\qbezier(2.65,0.9)(3.5,0.9)(3.5,0.9)
\put(0,2.7){$b_{i+1(2n)}$}
\put(1.3,2.8){\vector(0,1){0.01}}
\put(3.0,1.45){$a_{i-1(2n)}$}
\put(3.3,1.3){\vector(1,0){0.01}}
\put(-1.25,1.25){$a_i$}
\put(-0.8,1.35){\vector(0,1){0.01}}
\put(0.65,0.65){$b_i$}
\put(1.1,0.7){\vector(0,1){0.01}}
\qbezier[20](3,-2.5)(3,-2.25)(3,-1.5)
\put(3.01,-1.5){\vector(0,1){0.75}}
}\ \ \ \ \ \ \ \ \ \ \ \ \ \ \ \ 
\ \ \ \ \ \ \ \ \ .
\begin{array}{l}\\ \\ \\ \\ \\ \\ \\ \\ \\ \\ \\ \\  \end{array}
\end{eqnarray}
\setlength{\unitlength}{10pt}

Take the universal cyclic cover of $S^3\setminus U^{W_{2n}}$ by way of the surface
$\Sigma$ as in Subsection~{\rm\ref{uccs}. Build the covering space so that in 
the figure, if the projection of a path starting
in $X_i$ travels through the surface in the direction indicated then the
path is now in $X_{i+1}$. Choose the generating covering transformation as the
one which maps $X_{i}$ to $X_{i+1}$. 

We calculate the $\Zset[t,t^{-1}]$-module $H_1(\widetilde{S^3\setminus U^{W_{2n}}},\Zset)$.
$H_1(\widetilde{S^3\setminus U^{W_{2n}}},\Zset)$ has a generator for each component
of the surgery link. Thus it has a generator for each lift of the $a_i$ and the
$b_j$. Choose meridians $\alpha_i$ and $\beta_j$ having linking
number $+1$ with the associated components as oriented in the figure
to represent these generators. 

Relations are obtained by writing the lifts of the (0-framed) longitudes of
the $a_i$ and $b_j$ in terms of these meridians.
Up to translation there are $4n$ relations:

\begin{equation}
\begin{array}{lll}
\mbox{Around}\ \alpha_i& : & (t^{-1}-1)\beta_i + \beta_{i+1(2n)} = 0, \\
\mbox{Around}\ \beta_i& : & (t-1)\alpha_i + \alpha_{i-1(2n)} = 0,
\end{array}
\end{equation}
where $i$ runs from $1$ to $2n$ and in some cases has been written mod $2n$.

It follows from this that 
\begin{equation}
H_1(\widetilde{S^3\setminus U^{W_{2n}}},\Zset) =
\frac{\Zset[t,t^{-1}]}{(1-(1-t)^{2n})}\oplus \frac{\Zset[t,t^{-1}]}{(1-(1-t^{-1})^{2n})},
\end{equation}
and that the corresponding invariants are
\begin{eqnarray}
A_{U^{W_{2n}}}(t) & = & (1-(1-t)^{2n})(1-(1-t^{-1})^{2n}), \\
{C}_{U^{W_{2n}}}(h) & = & 1 - 2h^{2n} + \ldots ,
\end{eqnarray}
as is required for the proof of Equation~{\rm\ref{FC}}.

The form of $A$ above ($f(t)f(t^{-1})$) is that expected for a slice knot,
and in fact one can show that $U^{W_{2n}}$ is ribbon,
in fact a twisted form of Ng's wheels \cite{N}.

\section{Questions}

Notice that the result for diagrams whose graphs have 
negative Euler characteristic
is better than we could have expected. If $G$ is a graph clasper for
$U$ with underlying complete, primitive web diagram $D$ of degree $n$ with
negative Euler characteristic then
all that was required was that ${C}(U^{G})=1+o(h^{n+1})$. However
we have found that ${C}(U^G)=1$.

This prompts the following question, whose (positive) answer is part of a 
work in progress:
\begin{question}
With terms as in the above paragraph, is it true that
\begin{equation}
\KI(U^G) = 1\ \mbox{mod}\ \mbox{span}_{\Qset}\left\{D'\ \mbox{s.t.}\ \chi(D')\leq \chi(D)\right\}, 
\end{equation}
where $\chi(D)$ is the Euler characteristic of the dashed graph of $D$?
\end{question}

Further observe that the branched cyclic cover over a knot obtained by
claspering the unknot by a clasper graph, the dashed graph of whose underlying
diagram has negative Euler characteristic, is a $\Zset$-homology sphere. Does
its $n$-triviality reflect the graph of its diagram in some way?

\bibliographystyle{amsalpha}

\end{document}